# A Geometric Proof of Riemann Hypothesis

**Kaida Shi**

Department of Mathematics, Zhejiang Ocean University,

Zhoushan City, Zip.316004, Zhejiang Province, China

**Abstract** Beginning from the resolution of Riemann Zeta function $\zeta(s)$, using the inner product formula of infinite-dimensional vectors in the complex space, the author proved the world's problem Riemann hypothesis raised by German mathematician Riemannn in 1859.

MR(2000): 11M

## 1 Introduction

In 1859, German mathematician B. Riemann raised six very important hypothesizes about $\zeta(s)$ in his thesis entitled ***Ueber die Anzahl der Primzahlen unter einer gegebenen Große***[1]. Since his thesis was published for 30 years, five of the six hypothesizes have been solved satisfactorily. The remained one is the Super baffling problem—Riemann hypothesis. That is the non-trivial zeroes of the Riemann Zeta function $\zeta(s)$ are all lie on the straight line $\mathrm{Re}(s) = \frac{1}{2}$ within the complex plane $s$. The well-known Riemann Zeta function $\zeta(s)$ raised by Swiss mathematician Leonard Euler on 1730 to 1750. In 1749, he had proved the Riemann Zeta function $\zeta(s)$ satisfied another function equation. For $\mathrm{Re}(s) > 1$, the series

$$\frac{1}{1^s} + \frac{1}{2^s} + \frac{1}{3^s} + \cdots + \frac{1}{n^s} + \cdots$$

is convergent, it can be defined as $\zeta(s)$. **Although the definitional domain of the Riemann Zeta function $\zeta(s) = \sum_{n=1}^{\infty} \frac{1}{n^s}\ (s = \sigma + it)$ is $\mathrm{Re}(s) > 1$, but German mathematician Hecke used the** **analytic continuation method to continue this definitional domain to whole complex plane (except** $\text{Re}(s) = 1$**).**

In 1914, English mathematician G. H. Hardy declared that he had proved Riemann Zeta function $\zeta(s)$ equation has infinite non-trivial zeroes lie on the straight line $\text{Re}(s) = \frac{1}{2}$.

In fact, about the Riemann Zeta functional equality

$$\pi^{-\frac{s}{2}}\Gamma(\frac{s}{2})\zeta(s) = \pi^{-\frac{1-s}{2}}\Gamma(\frac{1-s}{2})\zeta(1-s), \tag{*}$$

because

$$\pi^{-\frac{\frac{1}{2}}{2}}\Gamma(\frac{\frac{1}{2}}{2})\zeta(\frac{1}{2}) = \pi^{-\frac{1-\frac{1}{2}}{2}}\Gamma(\frac{1-\frac{1}{2}}{2})\zeta(1-\frac{1}{2}),$$

therefore, **we have discovered that the solution of above equality is** $s = \frac{1}{2}$, namely $\text{Re}(s) = \frac{1}{2}$.

That is to say no matter how the value $t$ takes, $\sigma$ always equals to $\frac{1}{2}$. **This explains that the definitional domain of the Riemann Zeta function** $\zeta(s)$ **can be continued to at least** $\text{Re}(s) = \frac{1}{2}$.

## 2 The formal resolution of Riemann Zeta function $\zeta(s)$

Considering the geometric meaning of the Riemann Zeta function

$$\zeta(s) = \frac{1}{1^s} + \frac{1}{2^s} + \frac{1}{3^s} + \cdots + \frac{1}{n^s} + \cdots,$$

we will have the following figure:

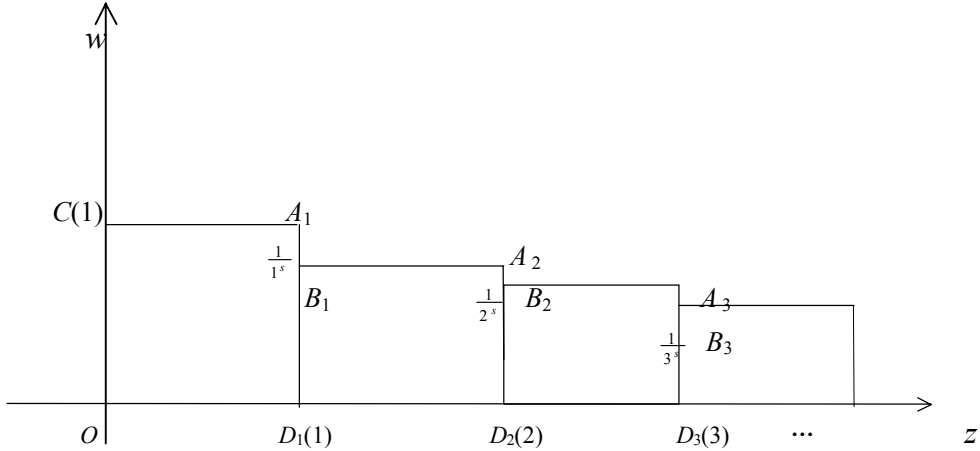

Fig. 1

In above figure, the areas of the rectangles $A_1D_1OC, A_2D_2D_1B_1, A_3D_3D_2B_2, \cdots$ are respectively:

$$\frac{1}{1^s}, \frac{1}{2^s}, \frac{1}{3^s}, \cdots,$$

therefore, the geometric meaning of the Riemann Zeta function $\zeta(s)$ is the sum of the areas of a series of rectangles within the complex space $s$.

Using the inner product formula between two infinite-dimensional vectors, the Riemann Zeta function $\zeta(s)$ equation

$$\zeta(s) = \frac{1}{1^s} + \frac{1}{2^s} + \frac{1}{3^s} + \cdots + \frac{1}{n^s} + \cdots = 0 \quad (s = \sigma + it)$$

can be resolved as

$$(\frac{1}{1^\sigma}, \frac{1}{2^\sigma}, \frac{1}{3^\sigma}, \cdots, \frac{1}{n^\sigma}, \cdots) \cdot (\frac{1}{1^{it}}, \frac{1}{2^{it}}, \frac{1}{3^{it}}, \cdots, \frac{1}{n^{it}}, \cdots) = 0.$$

So, we should have (**Please refer to the Appendix in the end**)

$$\sqrt{\frac{1}{1^{2\sigma}} + \frac{1}{2^{2\sigma}} + \frac{1}{3^{2\sigma}} + \cdots + \frac{1}{n^{2\sigma}} + \cdots} \cdot \sqrt{\frac{1}{1^{2it}} + \frac{1}{2^{2it}} + \frac{1}{3^{2it}} + \cdots \frac{1}{n^{2it}} + \cdots} \cos(\widehat{\vec{N}_1, \vec{N}_2}) = 0, \quad (1)$$

but in the complex space, if the inner product between two vectors equals to zero, then these two vectors are perpendicular, namely, $(\vec{N}_1, \vec{N}_2) = \dfrac{\pi}{2}$, therefore

$$\cos(\widehat{\vec{N}_1, \vec{N}_2}) = 0.$$

So long as the values of two radical expressions in the left side of (1) are finite (real or imaginary), then we think the expression (1) is tenable.

Because Riemann had proved that when $\mathrm{Re}(s) > 1$, the equation $\zeta(s) = 0$ has no solutions; when $\mathrm{Re}(s) < 0$, the Riemann Zeta function $\zeta(s)$ has only the trivial zeroes such as $s = -2,\ -4,\ -6,\ \cdots$, therefore, the non-trivial zeroes of the Riemann Zeta function $\zeta(s)$ will lie on the banded domain $0 \leq \sigma < 1$. **This explains that German mathematician Hecke used the analytic continuation method to continue the definitional domain of** t**he Riemann Zeta function** $\zeta(s)$ **is tenable.**

## 3 The relationship between the volume of the rotation solid and the area of its axis-cross section within the complex space

From the expression (1), we can know that when $\mathrm{Re}(s) = \sigma = \dfrac{1}{2}$, the series (harmonic series) within the first radical expression of left side is divergent. Because $\cos(\widehat{\vec{N}_1, \vec{N}_2}) = 0,$ therefore the situation of the complex series within the second radical expression of left side will change unable to research. Hence, we must transform the Riemann Zeta function $\zeta(s)$ equation. For this aim, let's derive the relationship between the volume of the rotation solid and the area of its axis-cross section within the complex space.

We call the cross section which pass through the axis $z$ and intersects the rotation solid as **axis-cross section**. According to the barycentric formula of the complex plane lamina:

$$\begin{cases} \xi = \dfrac{\int_a^b z[f(z)-g(z)]dz}{\int_a^b [f(z)-g(z)]dz}, \\ \eta = \dfrac{\frac{1}{2}\int_a^b [f^2(z)-g^2(z)]dz}{\int_a^b [f(z)-g(z)]dz}. \end{cases}$$

we have

$$2\pi\eta = \dfrac{\pi\int_a^b [f^2(z)-g^2(z)]dz}{\int_a^b [f(z)-g(z)]dz}. \tag{**}$$

The numerator of the fraction of right side of the formula (**) is the volume of the rotation solid, and the denominator is the area of the axis-cross section. The $\eta$ of left side of the formula (**) is the longitudinal coordinate of the barycenter of axis-cross section. Its geometric explanation is: taking the axis-cross section around the axis $z$ to rotate the angle of $2\pi$, we will obtain the volume of the rotation solid.

## 4 The proof of divergence of two concerned serieses

Now, let's prove the divergence of the series $\sum_{n=1}^{\infty} \cos(4t\ln n)$ and the series $\sum_{n=1}^{\infty} \sin(4t\ln n)$.

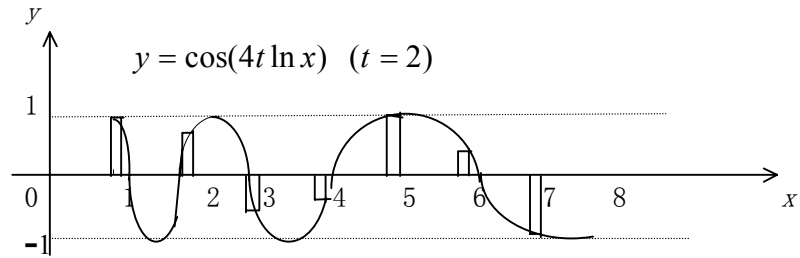

Fig. 2

Obviously, from Fig.2, we can find that the area of $c \cdot \sum_{n=1}^{\infty} \cos(4t\ln n)$ is contained by the area of $\int_{1}^{+\infty} \cos(4t \ln x)dx$.

Suppose that $c$ is a very small positive number, then the sum of the areas of a series of rectangles (which take $\cos(4t \ln n)(n = 1,2,3,\cdots)$ as **high** and take $c$ as **width**) is:

$$c \cdot \sum_{n=1}^{\infty} \cos(4t \ln n) = c \cdot \cos(4t \ln 1) + c \cdot \cos(4t \ln 2) + c \cdot \cos(4t \ln 3) + \cdots.$$

Because

$$\int_{1}^{+\infty} \cos(4t \ln x)dx = \lim_{\substack{w \to +\infty \\ t \to \infty}} \frac{1}{1+16t^2} x[\cos(4t \ln x) + 4t \sin(4t \ln x)]_{1}^{w}$$

$$= \lim_{\substack{w \to +\infty \\ t \to \infty}} \frac{1}{1+16t^2}[w(\cos(4t \ln w) + 4t \sin(4t \ln w)) - 1]$$

$$= \lim_{\substack{w \to +\infty \\ t \to \infty}} \frac{w}{\sqrt{1+16t^2}} \left( \frac{1}{\sqrt{1+16t^2}} \cos(4t \ln w) + \frac{4t}{\sqrt{1+16t^2}} \sin(4t \ln w) \right) - \lim_{t \to \infty} \frac{1}{1+16t^2}$$

Suppose that $\dfrac{1}{\sqrt{1+16t^2}} = \sin \delta$ and $\dfrac{4t}{\sqrt{1+16t^2}} = \cos \delta$, then above equals to

$$\lim_{\substack{w \to +\infty \\ t \to \infty}} \frac{w}{\sqrt{1+16t^2}} \sin(\delta + 4t \ln w) - \lim_{t \to \infty} \frac{1}{1+16t^2}$$

$$= \lim_{\substack{w \to +\infty \\ t \to \infty}} \frac{d}{\sqrt{\frac{1}{t^2}+16}} \sin(\delta + 4t \ln w)$$

$$= \frac{d}{4} \lim_{\substack{w \to +\infty \\ t \to \infty}} \sin(\delta + 4t \ln w).$$

Because the relationship between $w$ and $t$ is linear, therefore, we denote $w = dt$.

Obviously,

$$-\frac{d}{4} \le \frac{d}{4} \lim_{\substack{w \to +\infty \\ t \to \infty}} \sin(\delta + 4t \ln w) \le \frac{d}{4},$$

therefore

$$-\frac{d}{4} \le c \cdot \sum_{n=1}^{\infty} \cos(4t \ln n) \le \frac{d}{4},$$

namely,

$$-\frac{d}{4c} \le \sum_{n=1}^{\infty} \cos(4t \ln n) \le \frac{d}{4c}.$$

Thus, although the series $\sum_{n=1}^{\infty} \cos(4t \ln n)$ is divergent, but it is **oscillatory** and **bounded**.

Similarly, although the series $\sum_{n=1}^{\infty} \sin(4t \ln n)$ is divergent, but it is also **oscillatory** and **bounded**.

In a word, when $t = 0$, the series $\sum_{n=1}^{\infty} \cos(4t \ln n)$ is divergent (goes to **infinite**); when $t \ne 0$, it is **oscillatory** and **bounded**.

Thus, the divergence (**oscillatory** and **bounded**) of $\sum_{n=1}^{\infty} \cos(4t \ln n)$ has been proved.

Similarly, we can prove the divergence (**oscillatory** and **bounded**) of the series $\sum_{n=1}^{\infty} \sin(4t \ln n)$ $(t \ne 0)$.

.

## 5 Two theorems and their corollaries

Suppose that

$$\xi(s) = \frac{1}{2} s(s-1) \pi^{-\frac{1-s}{2}} \Gamma(\frac{1-s}{2}) \zeta(1-s),$$

then, we have

$$\xi(s) = \zeta(s).$$

According to the following:

**Theorem 1**[2] *The function $\xi(s)$ is the one class integer function, it has infinite zero points $\rho_n$, $0 \leq \mathrm{Re}(\rho_n) \leq 1$, the series $\sum_{n=1}^{\infty} |\rho_n|^{-1}$ is divergent. But for $\forall \varepsilon > 0$, the series $\sum_{n=1}^{\infty} |\rho_n|^{-1-\varepsilon}$ is convergent. The zeroes of the function $\xi(s)$ are non-trivial zeroes of the Riemann Zeta function.*

We can obtain:

**Corollary 1**[2] *The formula*

$$\xi(s) = e^{A+Bs} \prod_{n=1}^{\infty} \left(1 - \frac{s}{\rho_n}\right) e^{\frac{s}{\rho_n}}$$

*is tenable.*

**Corollary 2**[2] *The non-trivial zeroes of the Riemann Zeta function $\zeta(s)$ distribute symmetrically on the straight line $\mathrm{Re}\, s = \frac{1}{2}$ and $\mathrm{Im}\, s = 0$.*

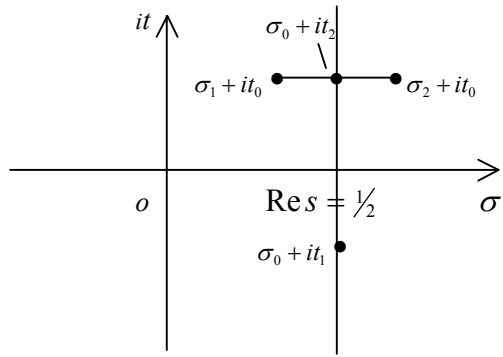

Fig. 3

According to the **symmetry** (about $\mathrm{Re}(s) = \frac{1}{2}$) of distribution of non-trivial zeroes of the Riemann Zeta function $\zeta(s)$ stated by the Corollary 2, we can consider the **uniqueness** of non-trivial zeroes of the Riemann Zeta function $\zeta(s)$ for both the imaginative coordinate $t_0$ and the real coordinate $\sigma_0$, we therefore suggest:

**Theorem 2 (The theorem about uniqueness and continuation of the non-trivial zero of the**

**Riemann Zeta function $\zeta(s)$)** *The imaginative coordinate $t_0$ of each non-trivial zero of the Riemann Zeta function $\zeta(s)$ doesn't correspond with two or over two real coordinates $\sigma_k (k = 1, 2, 3, \cdots)$; but the real coordinate $\sigma_0$ of each non-trivial zero of the Riemann Zeta function $\zeta(s)$ corresponds with infinite imaginative coordinates $t_k (k = 1, 2, 3, \cdots)$.*

**Proof** First, suppose that $t_0$ corresponds with two real coordinates $\sigma_1$ and $\sigma_2$, then, we have following equation group:

$$\begin{cases} \zeta(\sigma_1 + it_0) = \dfrac{1}{1^{\sigma_1+it_0}} + \dfrac{1}{2^{\sigma_1+it_0}} + \dfrac{1}{3^{\sigma_1+it_0}} + \cdots + \dfrac{1}{n^{\sigma_1+it_0}} + \cdots = 0, \\ \zeta(\sigma_2 + it_0) = \dfrac{1}{1^{\sigma_2+it_0}} + \dfrac{1}{2^{\sigma_2+it_0}} + \dfrac{1}{3^{\sigma_2+it_0}} + \cdots + \dfrac{1}{n^{\sigma_2+it_0}} + \cdots = 0. \end{cases}$$

Taking the first equation minus the second equation, we will obtain

$$\zeta(\sigma_1 + it_0) - \zeta(\sigma_2 + it_0) = \sum_{n=1}^{\infty} \left( \frac{1}{n^{\sigma_1+it_0}} - \frac{1}{n^{\sigma_2+it_0}} \right)$$

$$= \sum_{n=1}^{\infty} \frac{n^{\sigma_2} - n^{\sigma_1}}{n^{(\sigma_1+\sigma_2)+it_0}},$$

because $\sigma_1 \neq \sigma_2$, therefore, $n^{\sigma_2} - n^{\sigma_1} \neq 0$. On the other hand, because

$$n^{(\sigma_1+\sigma_2)+it_0} = n^{\sigma_1+\sigma_2} \cdot n^{it_0} = n^{\sigma_1+\sigma_2} \cdot e^{it_0 \ln n}$$

$$= n^{\sigma_1+\sigma_2} (\cos(t_0 \ln n) + i \sin(t_0 \ln n)) \neq 0,$$

therefore we have

$$\zeta(\sigma_1 + it_0) - \zeta(\sigma_2 + it_0) \neq 0$$

namely

$$\zeta(\sigma_1 + it_0) \neq \zeta(\sigma_2 + it_0).$$

This is contradictory with the fact that the first equation minus the second equation must equal to zero. Hence, the uniqueness of the real coordinate $\sigma$ of each non-trivial zero of the Riemann Zeta function $\zeta(s)$ has been proved.

Second, suppose that the real coordinate $\sigma_0$ of each non-trivial zero of the Riemann Zeta function $\zeta(s)$ corresponds with two imaginative coordinates $t_1$ and $t_2$, then, we have following equation group:

$$\begin{cases} \zeta(\sigma_0 + it_1) = \dfrac{1}{1^{\sigma_0+it_1}} + \dfrac{1}{2^{\sigma_0+it_1}} + \dfrac{1}{3^{\sigma_0+it_1}} + \cdots + \dfrac{1}{n^{\sigma_0+it_1}} + \cdots = 0, \\ \zeta(\sigma_0 + it_2) = \dfrac{1}{1^{\sigma_0+it_2}} + \dfrac{1}{2^{\sigma_0+it_2}} + \dfrac{1}{3^{\sigma_0+it_2}} + \cdots + \dfrac{1}{n^{\sigma_0+it_2}} + \cdots = 0. \end{cases}$$

Taking the first equation minus the second equation, we obtain

$$\zeta(\sigma_0 + it_1) - \zeta(\sigma_0 + it_2) = \sum_{n=1}^{\infty} \left( \frac{1}{n^{\sigma_0+it_1}} - \frac{1}{n^{\sigma_0+it_2}} \right)$$

$$= \sum_{n=1}^{\infty} \frac{n^{it_2} - n^{it_1}}{n^{\sigma_0+i(t_1+t_2)}} = \sum_{n=1}^{\infty} \frac{e^{it_2 \ln n} - e^{it_1 \ln n}}{n^{\sigma_0+i(t_1+t_2)}}$$

$$= \sum_{n=1}^{\infty} \frac{(\cos(t_2 \ln n) - \cos(t_1 \ln n)) + i(\sin(t_2 \ln n) - \sin(t_1 \ln n))}{n^{\sigma_0+i(t_1+t_2)}},$$

where

$$n^{\sigma_0+i(t_1+t_2)} = n^{\sigma_0} \cdot n^{i(t_1+t_2)} = n^{\sigma_0} \cdot e^{i(t_1+t_2)\ln n}$$
$$= n^{\sigma_0} \cdot (\cos(t_1 \ln n) + i\sin(t_1 \ln n)) \cdot (\cos(t_2 \ln n) + i\sin(t_2 \ln n))$$
$$\neq 0.$$

Enable above expression equals to zero, we must have

$$\begin{cases} \cos(t_2 \ln n) = \cos(t_1 \ln n), \\ \sin(t_2 \ln n) = \sin(t_1 \ln n) \end{cases} \quad (n = 1,\ 2,\ 3,\ \cdots)$$

so, we obtain $t_1 = t_2 + \dfrac{2k\pi}{\ln n}$ $(k = 1,\ 2,\ 3,\cdots)$. That is to say $t_1$ and $t_2$ can take any value, but

$$t_1 - t_2 = \frac{2k\pi}{\ln n} \quad (k = 1,\ 2,\ 3,\cdots).$$

Thus, theorem 2 has been proved.

**Corollary 3** *According to the* **solution** *of the equality* (*) *and the* **symmetry** *of Corollary 2 and*

*the **uniqueness** and the **continuation** of Theorem 2, we can say certainly that if $t \in (-\infty, +\infty)$, then all non-trivial zeroes of the Riemann Zeta function $\zeta(s)$ lie on straight line $\text{Re}(s) = \frac{1}{2}$ (except a certain point).*

## 6 The proof of Riemann hypothesis

Because we have the volume formula of the rotation solid formed by rotating the rectangle $A_n D_n D_{n-1} B_{n-1}$ around the axis $z$:

$$V_n = \int_{n-1}^{n} \frac{\pi}{n^{2s}} dz = \frac{\pi}{n^{2s}} z \Big|_{n-1}^{n} = \frac{\pi}{n^{2s}},$$

therefore, the sum of the volumes of a series of the cylinders formed by rotating a series of the rectangles $A_1 D_1 OC, A_2 D_2 D_1 B_1, A_3 D_3 D_2 B_2, \cdots$ around the axis $z$ is:

$$V = \sum_{n=1}^{\infty} V_n = \pi \left( \frac{1}{1^{2s}} + \frac{1}{2^{2s}} + \frac{1}{3^{2s}} + \cdots + \frac{1}{n^{2s}} + \cdots \right).$$

But the sum of the areas of a series of the rectangles $A_1 D_1 OC, A_2 D_2 D_1 B_1, A_3 D_3 D_2 B_2, \cdots$ is:

$$\frac{1}{1^s} + \frac{1}{2^s} + \frac{1}{3^s} + \cdots + \frac{1}{n^s} + \cdots.$$

By the formula (**) in Section 2, we have

$$\pi \left( \frac{1}{1^{2s}} + \frac{1}{2^{2s}} + \frac{1}{3^{2s}} + \cdots + \frac{1}{n^{2s}} + \cdots \right)$$

$$= 2\pi\eta \left( \frac{1}{1^s} + \frac{1}{2^s} + \frac{1}{3^s} + \cdots + \frac{1}{n^s} + \cdots \right). \tag{2}$$

Substituting the Riemann Zeta function $\zeta(s)$ equation

$$\zeta(s) = \frac{1}{1^s} + \frac{1}{2^s} + \frac{1}{3^s} + \cdots + \frac{1}{n^s} + \cdots = 0$$

into (2), we can obtain the transformed equation:

$$\frac{1}{1^{2s}} + \frac{1}{2^{2s}} + \frac{1}{3^{2s}} + \cdots + \frac{1}{n^{2s}} + \cdots = 0.$$

This can be explained geometrically as: **if the area of the axis-cross section equals to zero, then it corresponds to the volume of the rotation solid also equals to zero.**

Because $s = \sigma + it$, therefore

$$\frac{1}{1^{2(\sigma+it)}} + \frac{1}{2^{2(\sigma+it)}} + \frac{1}{3^{2(\sigma+it)}} + \cdots + \frac{1}{n^{2(\sigma+it)}} + \cdots = 0.$$

According to the inner product formula between two infinite-dimensional vectors, above expression can be written as:

$$(\frac{1}{1^{2\sigma}}, \frac{1}{2^{2\sigma}}, \frac{1}{3^{2\sigma}}, \cdots, \frac{1}{n^{2\sigma}}, \cdots) \cdot (\frac{1}{1^{2it}}, \frac{1}{2^{2it}}, \frac{1}{3^{2it}}, \cdots, \frac{1}{n^{2it}}, \cdots) = 0.$$

Denoting above two infinite-dimensional vectors as $\vec{N}_3$ and $\vec{N}_4$ respectively, we have

$$\sqrt{\frac{1}{1^{4\sigma}} + \frac{1}{2^{4\sigma}} + \frac{1}{3^{4\sigma}} + \cdots + \frac{1}{n^{4\sigma}} + \cdots} \cdot \sqrt{\frac{1}{1^{4it}} + \frac{1}{2^{4it}} + \frac{1}{3^{4it}} + \cdots \frac{1}{n^{4it}} + \cdots} \cdot \cos(\widehat{\vec{N}_3, \vec{N}_4}) = 0 \quad (3)$$

Now, let's prove Riemann hypothesis. If $t$ will take all values of a certain infinite interval, then we think that Riemann hypothesis has been proved.

From the Riemann Zeta functional equality

$$\pi^{-\frac{s}{2}} \Gamma(\frac{s}{2}) \zeta(s) = \pi^{-\frac{1-s}{2}} \Gamma(\frac{1-s}{2}) \zeta(1-s), \tag{*}$$

we have known its solution is $\text{Re}(s) = \frac{1}{2}$. Substituting it into the expression (3), we will have:

$$\sqrt{\frac{1}{1^2} + \frac{1}{2^2} + \frac{1}{3^2} + \cdots + \frac{1}{n^2} + \cdots} \cdot \sqrt{\frac{1}{1^{4it}} + \frac{1}{2^{4it}} + \frac{1}{3^{4it}} + \cdots + \frac{1}{n^{4it}} + \cdots} \cdot \cos(\widehat{\vec{N}_3, \vec{N}_4}) = 0$$

namely,

$$\frac{\pi}{\sqrt{6}}\sqrt{\frac{1}{1^{4it}}+\frac{1}{2^{4it}}+\frac{1}{3^{4it}}+\cdots+\frac{1}{n^{4it}}+\cdots}\cdot\cos(\widehat{\vec{N}_3,\vec{N}_4})=0.$$

When $(\widehat{\vec{N}_3,\vec{N}_4})=\dfrac{\pi}{2}$, we have

$$\frac{\pi}{\sqrt{6}}\sqrt{\frac{1}{1^{4it}}+\frac{1}{2^{4it}}+\frac{1}{3^{4it}}+\cdots+\frac{1}{n^{4it}}+\cdots}\cdot 0=0.$$

In above expression, when $t=0$, although the series

$$\sum_{n=1}^{\infty}\sin(4t\ln n)=\sum_{n=1}^{\infty}\sin(4\cdot 0\cdot\ln n)$$
$$=0+0+0+\cdots+0+\cdots=0,$$

but because the series

$$\sum_{n=1}^{\infty}\cos(4t\ln n)=\sum_{n=1}^{\infty}\cos(4\cdot 0\cdot\ln n)$$
$$=1+1+1+\cdots+1+\cdots=\infty,$$

therefore, the complex series

$$\frac{1}{1^{4it}}+\frac{1}{2^{4it}}+\frac{1}{3^{4it}}+\cdots+\frac{1}{n^{4it}}+\cdots$$
$$=\sum_{n=1}^{\infty}[\cos(4t\ln n)-i\sin(4t\ln n)]$$
$$=(1+1+1+\cdots+1+\cdots)-$$
$$-i\cdot(0+0+0+\cdots+0+\cdots)=\infty.$$

Because "infinite" multiplied by zero doesn't equal to zero, therefore

$$\frac{\pi}{\sqrt{6}}\sqrt{\frac{1}{1^{4it}}+\frac{1}{2^{4it}}+\frac{1}{3^{4it}}+\cdots+\frac{1}{n^{4it}}+\cdots}\cdot 0$$
$$=\frac{\pi}{\sqrt{6}}\cdot\sqrt{\sum_{n=1}^{\infty}[\cos(4t\ln n)-i\sin(4t\ln n)]}\cdot 0$$
$$=\frac{\pi}{\sqrt{6}}\cdot\sqrt{\infty}\cdot 0\neq 0.$$

So, $\text{Re}(s) = \dfrac{1}{2}$ $(t = 0)$ is not a non-trivial zero of the Riemann Zeta function $\zeta(s)$, it must be removed.

When $t \neq 0$, although the series $\sum\limits_{n=1}^{\infty} \cos(4t \ln n)$ and series $\sum\limits_{n=1}^{\infty} \sin(4t \ln n)$ are all divergent (see Section 3), but they are all **oscillatory** and **bounded**, therefore we can know the complex series

$$\dfrac{1}{1^{4it}} + \dfrac{1}{2^{4it}} + \dfrac{1}{3^{4it}} + \cdots + \dfrac{1}{n^{4it}} + \cdots$$
$$= \sum_{n=1}^{\infty}[\cos(4t \ln n) - i \sin(4t \ln n)]$$

is also **oscillatory** and **bounded**.

From following figure:

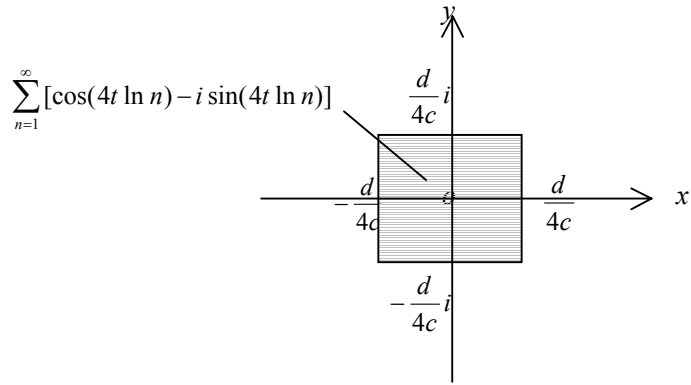

Fig. 4

we can know the change scope of the complex series $\sum\limits_{n=1}^{\infty}[\cos(4t \ln n) - i \sin(4t \ln n)]$ is the square.

Obviously, the finite complex numbers multiplied by zero equal to zero, therefore we have

$$\frac{\pi}{\sqrt{6}}\sqrt{\frac{1}{1^{4it}}+\frac{1}{2^{4it}}+\frac{1}{3^{4it}}+\cdots+\frac{1}{n^{4it}}+\cdots}\cdot 0$$

$$=\frac{\pi}{\sqrt{6}}\cdot\sqrt{\sum_{n=1}^{\infty}[\cos(4t\ln n)-i\sin(4t\ln n)]}\cdot 0$$

$$=0, \qquad\qquad (t\neq 0)$$

By the expression (1), considering the situation of $\Gamma(\frac{s}{2})$ in $\zeta(s)$, we have known that when $s = 0$, $\zeta(s) \neq 0$, therefore, $s = 0$ is not a non-trivial zero of the Riemann Zeta function $\zeta(s)$. Because $s = \sigma + it$, therefore $t = 0$ is also not a non-trivial zero of the Riemann Zeta function $\zeta(s)$. This is identical with the situation of the sufficient and essential condition of divergence of the complex series, therefore, we obtain

$$t \in (-\infty,\ 0) \cup (0,\ +\infty).$$

**This explains that $t$ takes all values of the straight line $\text{Re}(s) = \frac{1}{2}$ except $t = 0$.**

Obviously, this result is identical with the corollaries of theorem 1 and theorem 2 within Section 4.

## 7 The progress in the process for proving Riemann hypothesis

Riemann has proved that Zeta function satisfied the function equation:

$$\zeta(1-s) = \zeta(s)\Gamma(\frac{s}{2})\pi^{-\frac{s}{2}},$$

in which, when $s$ is the positive number, we have

$$\Gamma(s) = \int_0^{\infty} x^{s-1}e^{-x}dx.$$

Especially, when $s$ is the positive integral, we have

$$\Gamma(s+1) = s!,$$

when *s* is arbitrary number (include negative number), we have

$$\Gamma(s) = \lim_{n \to \infty} \frac{n! n^{s-1}}{s(s+1)(s+2)\cdots(s+n-1)}.$$

In 1942, Norway mathematician A. Selberg proved

$$N_0(T) \geq CT \log T, \quad (C > 0).$$

In fact, Selberg proved

$$N_0(T) \geq CN(T), \quad (C > 0, \ C \in R).$$

but, Selberg obtained the constant *C* has only 0.01.

In 1974, mathematician N. Levinson at Massachusetts Institute of Technology of USA proved successfully that for full great real number *T*,

$$N_0(T) \geq \frac{1}{3} N(T).$$

that is to say Riemann Zeta function $\zeta(s)$ has at least a third non-trivial zeroes lie on the straight line $\text{Re}(s) = \frac{1}{2}$.

In 1980, Chinese mathematician Shituo Lou after improving slightly Levinson's result and proved:

$$N_0(T) > 0.35 N(T).$$

These are most new general survey about the research of Riemann hypothesis[3~8].

In 1968, by using great computer, three mathematicians at University of Wisconsin of USA proved first 3,000,000 non-trivial zeroes of Riemann Zeta function $\zeta(s)$ lie on the straight line $\text{Re}(s) = \frac{1}{2}$. After that, even one man using great computer checked the first 1,500,000,000 non-trivial zeroes of Riemann Zeta function $\zeta(s)$ lie on the straight line $\text{Re}(s) = \frac{1}{2}$[8]. But the ability of the computer is limited, it can't achieve infinite, therefore, that can't to say whether Riemann hypothesis about the non-trivial zeroes of Zeta function $\zeta(s)$ is correct or not, it only provides the basis for supporting the hypothesis.

**References**

1. Lou S. T., Wu D. H., *Riemann hypothesis*, Shenyang: Liaoning Education Press, 1987. pp.152-154.

2. Chen J. R., On the representation of a large even integer as the sum of a prime and the product of at most two primes, *Science in China*, 16 (1973), No.2, pp. 111-128..

3. Pan C. D., Pan C. B., *Goldbach hypothesis*, Beijing: Science Press, 1981. pp.1-18; pp.119-147.

4. Hua L. G., *A Guide to the Number Theory*, Beijing: Science Press, 1957. pp.234-263.

5. Chen J. R., Shao P. C., *Goldbach hypothesis*, Shenyang: Liaoning Education Press, 1987. pp.77-122; pp.171-205.

6. Chen J. R., *The Selected Papers of Chen Jingrun*, Nanchang: Jiangxi Education Press, 1998. pp.145-172.

7. Lehman R. S., On the difference $\pi(x)$-li$x$, *Acta Arith.*, 11(1966). pp.397～410.

8. Hardy, G. H., Littlewood, J. E., Some problems of "patitio numerorum" III: On the expression of a number as a sum of primes, *Acta. Math.*, 44 (1923). pp.1-70.

9. Hardy, G. H., Ramanujan, S., Asymptotic formula in combinatory analysis, *Proc. London Math. Soc.*, (2) 17 (1918). pp. 75-115.

10. Riemann, B., Ueber die Anzahl der Primzahlen unter einer gegebenen Große, Ges. Math. Werke und Wissenschaftlicher Nachlaß , 2, Aufl, 1859, pp 145-155.

11. E. C. Titchmarsh, *The Theory of the Riemann Zeta Function*, Oxford University Press, New York, 1951.

12. Morris Kline, *Mathematical Thought from Anoient to Modern Times*，Oxford University Press, New York, 1972.

13. A. Selberg, *The Zeta and the Riemann Hypothesis,* Skandinaviske Matematiker Kongres, 10 (1946).

# **Appendix**

## **The inner product formula between two vectors in real space can be extended formally to complex space**

**Dear Mr. Referee,**

**Thank you for your review.**
**First, please observe following examples.**

**Example 1** $A\{1+i,\ 3\},\ B\{-i,\ 2i\},\ C\{1,\ -i\}$.

$\overrightarrow{AB} = \{-1-2i,\ -3+2i\},\ \overrightarrow{AC} = \{-i,\ -3-i\},\ \overrightarrow{BC} = \{1+i,\ -3i\}$.

$\left|\overrightarrow{AB}\right| = \sqrt{(-1-2i)^2 + (-3+2i)^2} = \sqrt{2-8i}$,

$\left|\overrightarrow{AC}\right| = \sqrt{(-i)^2 + (-3-i)^2} = \sqrt{7+6i}$,

$\left|\overrightarrow{BC}\right| = \sqrt{(1+i)^2 + (-3i)^2} = \sqrt{-9+2i}$.

According to **Cosine theorem**

$$\left|\overrightarrow{AB}\right|^2 + \left|\overrightarrow{AC}\right|^2 - \left|\overrightarrow{BC}\right|^2 = 2\left|\overrightarrow{AB}\right|\left|\overrightarrow{AC}\right|\cos(\overrightarrow{AB},\overrightarrow{AC}),$$

we have

$$\frac{\left|\overrightarrow{AB}\right|^2 + \left|\overrightarrow{AC}\right|^2 - \left|\overrightarrow{BC}\right|^2}{2} = \frac{(2-8i)+(7+6i)-(-9+2i)}{2} = 9-2i.$$

On the other hand, we have

$$\overrightarrow{AB}\cdot\overrightarrow{AC} = \left|\overrightarrow{AB}\right|\left|\overrightarrow{AC}\right|\cos(\overrightarrow{AB},\overrightarrow{AC}) = (-1-2i)(-i)+(-3+2i)(-3-i) = 9-2i.$$

According to the **Cosine theorem**

$$\left|\overrightarrow{AC}\right|^2 + \left|\overrightarrow{BC}\right|^2 - \left|\overrightarrow{AB}\right|^2 = 2\left|\overrightarrow{AC}\right|\left|\overrightarrow{BC}\right|\cos(\overrightarrow{AC},\overrightarrow{BC})$$

$$\frac{\left|\overrightarrow{AC}\right|^2 + \left|\overrightarrow{BC}\right|^2 - \left|\overrightarrow{AB}\right|^2}{2} = \frac{(7+6i)+(-9+2i)-(2-8i)}{2} = -2+8i.$$

On the other hand, we have

$$\overrightarrow{AC}\cdot\overrightarrow{BC} = \left|\overrightarrow{AB}\right|\left|\overrightarrow{AC}\right|\cos(\overrightarrow{AB},\overrightarrow{AC}) = (-i)(1+i)+(-3-i)(-3i) = -2+8i.$$

$$S_{\Delta ABC} = \frac{1}{2}\left|\overrightarrow{AB}\right|\left|\overrightarrow{AC}\right|\sin\theta_1 = \frac{1}{2}\sqrt{2-8i}\cdot\sqrt{7+6i}\cdot\frac{\sqrt{-15-8i}}{\sqrt{62-44i}} = \frac{1}{2}\sqrt{-15-8i},$$

$$S_{\Delta ACB} = \frac{1}{2}\left|\overrightarrow{AC}\right|\left|\overrightarrow{BC}\right|\sin\theta_{21} = \frac{1}{2}\sqrt{7+6i}\cdot\sqrt{-9+2i}\cdot\frac{\sqrt{-15-8i}}{\sqrt{-75-40i}} = \frac{1}{2}\sqrt{-15-8i}.$$

**Example 2** $A\{1+i,\ 1-i,\ 2i\}$, $B\{1-i,\ 1+i,\ -2i\}$, $C\{1,\ 0,\ i\}$.

$\overrightarrow{AB} = \{-2i,\ 2i,\ -4i\}$, $\overrightarrow{AC} = \{-i,\ -1+i,\ -i\}$, $\overrightarrow{BC} = \{i,\ -1-i,\ 3i\}$.

$\left|\overrightarrow{AB}\right| = \sqrt{(-2i)^2 + (2i)^2 + (-4i)^2} = \sqrt{-24},$

$\left|\overrightarrow{AC}\right| = \sqrt{(-i)^2 + (-1+i)^2 + (-i)^2} = \sqrt{-2-2i},$

$\left|\overrightarrow{BC}\right| = \sqrt{(i)^2 + (-1-i)^2 + (3i)^2} = \sqrt{-10+2i}.$

According to the **Cosine theorem**

$$\left|\overrightarrow{AB}\right|^2 + \left|\overrightarrow{AC}\right|^2 - \left|\overrightarrow{BC}\right|^2 = 2\left|\overrightarrow{AB}\right|\left|\overrightarrow{AC}\right|\cos(\overrightarrow{AB},\overrightarrow{AC}),$$

we have

$$\frac{\left|\overrightarrow{AB}\right|^2 + \left|\overrightarrow{AC}\right|^2 - \left|\overrightarrow{BC}\right|^2}{2} = \frac{(-24) + (-2-2i) - (-10+2i)}{2} = -8 - 2i.$$

On the other hand, we have

$$\overrightarrow{AB}\cdot\overrightarrow{AC} = \left|\overrightarrow{AB}\right|\left|\overrightarrow{AC}\right|\cos(\overrightarrow{AB},\overrightarrow{AC}) = (-2i)(-i) + (2i)(-1+i) + (-4i)(-i) = -8 - 2i.$$

According to the **Cosine theorem**

$$\left|\overrightarrow{AC}\right|^2 + \left|\overrightarrow{BC}\right|^2 - \left|\overrightarrow{AC}\right|^2 = 2\left|\overrightarrow{AC}\right|\left|\overrightarrow{BC}\right|\cos(\overrightarrow{AC},\overrightarrow{BC}),$$

we have

$$\frac{\left|\overrightarrow{AC}\right|^2 + \left|\overrightarrow{BC}\right|^2 - \left|\overrightarrow{AB}\right|^2}{2} = \frac{(-2-2i) + (-10+2i) - (-24)}{2} = 6.$$

On the other hand, we have

$$\overrightarrow{AC}\cdot\overrightarrow{BC} = \left|\overrightarrow{AC}\right|\left|\overrightarrow{BC}\right|\cos(\overrightarrow{AC},\overrightarrow{BC}) = (-i)(i) + (-1+i)(-1-i) + (-i)(3i) = 6.$$

$$S_{\Delta ABC} = \frac{1}{2}\left|\overrightarrow{AB}\right|\left|\overrightarrow{AC}\right|\sin\theta_1 = \frac{1}{2}\sqrt{-24}\cdot\sqrt{-2-2i}\cdot\frac{\sqrt{-12+16i}}{\sqrt{48+48i}} = \sqrt{-3+4i},$$

$$S_{\Delta ACB} = \frac{1}{2}\left|\overrightarrow{AC}\right|\left|\overrightarrow{BC}\right|\sin\theta_{21} = \frac{1}{2}\sqrt{-2-2i}\cdot\sqrt{-10+2i}\cdot\frac{\sqrt{-12+16i}}{\sqrt{24+16i}} = \sqrt{-3+4i}.$$

**Example 3** $A\{8i,\ 14,\ 8-i,\ 1\},\ B\{6,\ 15i,\ 17,\ -8\},\ C\{3-i,\ 10+7i,\ 11,\ 3i\}$.

$\overrightarrow{AB} = \{6-8i,\ -14+15i,\ 9+i,\ -9\},\ \overrightarrow{AC} = \{3-9i,\ -4+7i,\ 3+i,\ -1+3i\},$
$\overrightarrow{BC} = \{-3-i,\ 10-8i,\ -6,\ 8+3i\}.$

$\left|\overrightarrow{AB}\right| = \sqrt{(6-8i)^2 + (-14+5i)^2 + (9+i)^2 + (-9)^2} = \sqrt{104 - 498i},$

$\left|\overrightarrow{AC}\right| = \sqrt{(3-9i)^2 + (-4+7i)^2 + (3+i)^2 + (-1+3i)^2} = \sqrt{-105-110i},$

$\left|\overrightarrow{BC}\right| = \sqrt{(-3-i)^2 + (10-8i)^2 + (-6)^2 + (8+3i)^2} = \sqrt{135 - 106i}.$

According to the **Cosine theorem**

$$\left|\overrightarrow{AB}\right|^2 + \left|\overrightarrow{AC}\right|^2 - \left|\overrightarrow{BC}\right|^2 = 2\left|\overrightarrow{AB}\right|\left|\overrightarrow{AC}\right|\cos(\overrightarrow{AB},\overrightarrow{AC}),$$

we have

$$\frac{\left|\overrightarrow{AB}\right|^2 + \left|\overrightarrow{AC}\right|^2 - \left|\overrightarrow{BC}\right|^2}{2} = \frac{(104-498i)+(-105-110i)-(-135-106i)}{2} = -68 - 251i.$$

On the other hand, we have

$$\overrightarrow{AB}\cdot\overrightarrow{AC} = \left|\overrightarrow{AB}\right|\left|\overrightarrow{AC}\right|\cos(\overrightarrow{AB},\overrightarrow{AC}) =$$
$$= (6-8i)(3-9i) + (-14+15i)(-4+7i) + (9+i)(3+i) + (-9)(-1+3i) = -68 - 251i.$$

According to the **Cosine theorem**

$$\left|\overrightarrow{AC}\right|^2 + \left|\overrightarrow{BC}\right|^2 - \left|\overrightarrow{AC}\right|^2 = 2\left|\overrightarrow{AC}\right|\left|\overrightarrow{BC}\right|\cos(\overrightarrow{AC},\overrightarrow{BC}),$$

we have

$$\frac{\left|\overrightarrow{AC}\right|^2 + \left|\overrightarrow{BC}\right|^2 - \left|\overrightarrow{AB}\right|^2}{2} = \frac{(-105-110i)+(135-106i)-(104-498i)}{2} = -37 + 141i.$$

On the other hand, we hand

$$\overrightarrow{AC} \cdot \overrightarrow{BC} = \left|\overrightarrow{AC}\right|\left|\overrightarrow{BC}\right|\cos(\overrightarrow{AC},\overrightarrow{BC}) =$$
$$= (3-9i)(-3-i) + (-4+7i)(10-8i) + (-6)(3+i) + (-1+3i)(8+3i) = -37 + 141i.$$

$$S_{\Delta ABC} = \frac{1}{2}\left|\overrightarrow{AB}\right|\left|\overrightarrow{AC}\right|\sin\theta_1 = \frac{1}{2}\sqrt{104-498i} \cdot \sqrt{-105-110i} \cdot \frac{\sqrt{-7323+6714i}}{\sqrt{-65700+40850i}}$$
$$= \frac{1}{2}\sqrt{-7323+6714i},$$

$$S_{\Delta ACB} = \frac{1}{2}\left|\overrightarrow{AC}\right|\left|\overrightarrow{BC}\right|\sin\theta_{21} = \frac{1}{2}\sqrt{-105-110i} \cdot \sqrt{135-106i} \cdot \frac{\sqrt{-7323+6714i}}{\sqrt{-25835-3720i}}$$
$$= \frac{1}{2}\sqrt{-7323+6714i}.$$